\documentclass[12pt]{amsart}
\usepackage{amsmath,amsthm,amsfonts,amssymb}
\usepackage[all]{xy}
\usepackage{amssymb, amsthm, amsmath}
\usepackage[english]{babel}
\usepackage{amsfonts}
\usepackage{tikz}

\DeclareMathOperator{\id}{id}
\DeclareMathOperator{\Hom}{Hom}

\newcommand{\ra}{\rightarrow}

\newcommand{\C}{\mathbb C}
\newcommand{\Z}{\mathbb Z}
\newcommand{\ot}{\otimes}
\newcommand{\mtc}{\mathcal}
\newcommand{\lam}{\lambda}

\newcommand{\al}{\alpha}
\newcommand{\eps}{\epsilon}
\newcommand{\bn}{\begin}

\newcommand{\D}{\Delta}

\newcommand{\rh}{\rightharpoonup}
\newcommand{\lh}{\leftharpoonup}
\numberwithin{equation}{section}
\newtheorem{lemma}[equation]{Lemma}

\newtheorem{defn}[equation]{Definition}
\newtheorem{cor}[equation]{Corollary}
\newtheorem{rem}[equation]{Remark}
\newcommand{\dw}{\downarrow}
\newcommand{\uw}{\uparrow}

\newcommand{\bl}{\begin{lemma}
  }
\newcommand{\nc}{\newcommand}
\nc{\el}{\end{lemma}}

\newcommand{\ch}{\chi}
\newcommand{\mtr}{\mathrm}
\nc{\bwt}{\bowtie}
\newcommand{\ncm}{\newcommand}\newcommand{\gm}{\gamma}
\numberwithin{equation}{section}
\newcommand{\et}{\end{thm}}\newcommand{\bt}{\bn{thm}}
\newcommand{\ep}{\end{prop}}\newcommand{\bp}{\bn{prop}}
\newcommand{\beqarn}{\begin{eqnarray*}}
\newcommand{\eeqarn}{\end{eqnarray*}}
\newcommand{\beqn}{\bn{equation*}}
\newcommand{\eeqn}{\end{equation*}}
\newcommand{\bpf}{\bn{proof}}
\newcommand{\epf}{\end{proof}}
\ncm{\cX}{\mtc{X}}
\ncm{\wt}{\widetilde}
\ncm{\sg}{\sigma}\ncm{\Rep}{\mathrm{Rep}}

\ncm{\X}{\mathcal{X}}
\ncm{\cA}{\mathcal{A}}
\newcommand{\lb}{\label}

\numberwithin{equation}{section}

\numberwithin{equation}{section}

\ncm{\np}{\newpage}
\ncm{\ebl}{\end{thebibliography}}
\ncm{\bbl}{\begin{thebibliography}}
\ncm{\chd}{_{ _{\ch}}}
\ncm{\ald}{_{ _{\al}}}
\ncm{\cP}{\mathcal{P}}
\ncm{\ei}{e_i}
\ncm{\eij}{e_{i,\;j}}
\ncm{\bne}{\begin{enumerate}}
\ncm{\ene}{\end{enumerate}}\ncm{\bdef}{\begin{defn}}
\ncm{\edf}{\end{defn}}
\ncm{\stab}{\mtr{Stab}}
\ncm{\bc}{\begin{cor}}

\ncm{\ec}{\end{cor}}
\ncm{\er}{\end{rem}}
\ncm{\br}{\begin{rem}}

\ncm{\bd}{\begin{document}}
\ncm{\ed}{\end{document}}
\ncm{\beq}{\begin{equation}}

\ncm{\eeq}{\end{equation}}
\ncm{\cm}{\mathcal{M}}
\ncm{\rep}{\mtr{Rep}}
\ncm{\btw}{\bowtie}
\ncm{\cd}{\mtc{D}}
\ncm{\cop}{\mtr{cop}}

\ncm{\bea}{\begin{eqnarray}}
\ncm{\eea}{\end{eqnarray}}
\ncm{\beanon}{\begin{eqnarray*}}
\ncm{\eeanon}{\end{eqnarray*}}\ncm{\ek}{\eps|_K}\ncm{\diez}{\#}

\ncm{\cC}{\mtc{C}}

\ncm{\cc}{\mtc{C}}

\ncm{\HKer}{\mtr{HKer}}
\ncm{\LKER}{\mtr{LKER}}
\ncm{\aad}{\mtr{ad}}
\ncm{\Dr}{\mtr{D}}
\ncm{\cD}{\mathcal{D}}
\ncm{\G}{\mathcal{G}}
\ncm{\Dc}{\mtc{D}}
\ncm{\E}{\mtc{E}}
\ncm{\fp}{\mtr{FPdim}}
\ncm{\Vc}{\mtr{Vec}}
\ncm{\cK}{\mtc{K}}
\ncm{\cM}{\mtc{M}}
\ncm{\cE}{\mtc{E}}
\ncm{\cS}{\mtc{S}}
\ncm{\End}{\mtr{End}}
\ncm{\subb}{\subsection{}\;\;}
\ncm{\subsubb}{\subsubsection{}\;\;\md}
\ncm{\md}{\medbreak}\ncm{\xra}{\xrightarrow}
\ncm{\hsa}{Hopf subalgebra of }
\ncm{\ses}{semisimple}
\ncm{\x}{$}
\ncm{\mi}{\mtr{I}}\ncm{\mrc}{\mtr{C}}
\ncm{\cZ}{\mtc{Z}}\ncm{\Gm}{\Gamma}
\ncm{\cb}{\mtc{B}}\ncm{\ca}{\mtc{A}}\ncm{\op}{\mtr{op}}
\ncm{\irr}{\Irr}\ncm{\Irr}{\mathrm{Irr}}
\ncm{\co}{\mtc{O}}%\ncm{\co}{\mtc{O}}
\ncm{\cg}{{\mtr{K}_0}}\ncm{\ci}{\mtc{I}}
\ncm{\blue}{\textcolor[rgb]{.00, .00, 1.00}}
\ncm{\red}{\textcolor[rgb]{1.00, .00, .00}}
\ncm{\green}{\textcolor[rgb]{.00, 1.00, .00}}
\numberwithin{equation}{section}
\title[Semisimple Hopf algebras]
{On the irreducible representations of generalized quantum doubles}%double crossed products via Clifford theory}
\author{Sebastian  Burciu}
\address{Inst.\ of Math.\ ``Simion Stoilow" of the Romanian Academy
P.O. Box 1-764, RO-014700, Bucharest, Romania}\address{and}
\address{University of Bucharest, Faculty of Mathematics and Computer Science, Algebra and Number Theory Research Center, 14 Academiei St., Bucharest, Romania }\email{sebastian.burciu@imar.ro} 
\thanks{This work was partially supported by a grant of the Romanian National Authority for Scientific Research, CNCS-UEFSCDI, project number PN II - ID - PCE-2011-3-0039}%{This research was partially supported by CNCSIS PN II RU-PD grant "Representations of semisimple Hopf algebras and fusion categories" PD-168 no: 14/28.07.2010}
\begin{document}
\subjclass[2000]{Primary 16W30, 18D10}
%\keywords{Double crossed products; Drinfeld doubles; Fusion categories; Clifford theory}
\keywords{Generalized quantum doubles; Fusion subcategories; Clifford Theory;  Grothendieck rings}

\maketitle
%\contents
\begin{abstract}
A description of all the irreducible representations of generalized quantum doubles associated to skew pairings of semisimple Hopf algebras is given. In particular a description of the irreducible representations of semisimple Drinfeld doubles is obtained.  It is shown that the Grothendieck ring of these generalized quantum doubles have a structure similar to the rings that arise from Green functors. In order to do this we give a formula for the tensor product of any two such irreducible representations.
\end{abstract}
%\tableofcontents\newpage
\section{Introduction and main results}

In what follows $k$ is an algebraically closed field of characteristic zero. All the algebras and fusion categories considered here are over $k$.

To any skew pairing $\lam:U \ot H \ra k$ one can associate a twisted Hopf algebra $D_{\lam}(U,H)$ of $U \otimes H$ called the generalized quantum double of $U$ and $H$ (see \cite{Maj}). If $H$ is finite dimensional, $U=H^{*\;\cop}$, and $\lam$ is the usual evaluation map then one obtains the usual quantum double $D(H)$ introduced by Drinfeld.

In this paper we describe all the irreducible representations of the Hopf algebras $D_{\lam}(U,H)$ when $U$ and $H$ are semisimple Hopf algebras and $\lam$ induces a surjective map $\lam(u,\;-):U \ra H^{*\;\cop}$. To this end, we apply Clifford's correspondence for semisimple Hopf algebras that was developed in \cite{coc}.

\md
The irreducible representations of the Drinfeld double $D(G)$ of a finite group are parameterized by pairs $(g,M)$ where $g$ runs through representatives of conjugacy classes of $G$, and $M$ is an irreducible $\mtr{C}_G(g)$-representation. To the pair $(g,M)$ it corresponds the induced representations $kG\ot_{k\mrc_G(g)}M$ (see \cite{L87} or \cite{W96}.)

Our main result is a description of irreducible representations of $D_{\lam}(U,H)$ that generalizes the above description for finite groups $G$. More precisely, we prove the following.

\bt\lb{moddescr} Suppose that $f :U \ra H^{*\;\cop}$ is a surjective Hopf algebra isomorphism. Then all the irreducible $D_{\lam_f}(U, H)$-representations are of the form
$S_{g,M}:=H\ot_{L(g)}M$ for some element $g \in G$ and some irreducible representation $M \in \ci_g$.
\et
\md
Here the group $G$ is the universal grading group of $\rep(U)$. Therefore $G$ is defined by $K(U)=kG^*$ where $K(U)$ is the largest central Hopf subalgebra of $U$. For the definition of the set $\ci_g$ see Equation \ref{i}. A definition of the Hopf subalgebra $L(g)$ of $H$ is given in \ref{defl} and it can be regarded as a generalization of the group algebra of the centralizer of a group element. Moreover, in Corollary \ref{bij} we also give a parametrization $(g, M)\mapsto S_{g,M}$ of all irreducible representations of $D_{\lam_f}(U,H)$ similar to the case of the quantum double $D(G)$.
\md
A similar description of the simple objects of the braided center of a graded fusion category was given recently in \cite{gnn}. It was shown in Proposition 3.9 of \cite{gnn} that in this situation the simple objects of $\cZ(\cc)$ are in bijection with the pairs formed by a representative of a conjugacy class of the universal grading group $U(\cc)$ and some simple equivariant objects of the corresponding grading components. Connections between the results obtained in \cite{gnn} and our results are provided in Section \ref{cocentral}.
\md
In Section \ref{catint} we give a categorical description of the category of representations $\rep(D_{\lam_f}(U,H))$ for a surjective morphism $f$ of Hopf algebras. We show that $\rep(D_{\lam_f}(U,H))$ is the relative center of the image functor $f_{*}$ induced by $f$. The notion of relative center of a tensor functor  appears in \cite{BN2} and is a particular case of the notion of relative center of a bimodule category introduced in \cite{gnn}.

\md
Our second main result is a formula for the tensor product $S_{g,M}\ot S_{h,N}$ of two such irreducible $D_{\lam_f}(U,H)$-representations. This is given in Section \ref{fus} where the following is proven:

\bt\lb{fusion} \lb{tp}Suppose that $f:U \ra H^{*\;\cop}$ is a surjective morphism of Hopf algebras.
Let $M\in \ci_g$ and $N \in \ci_h$ as above. Then the tensor product of two irreducible $D_f(U,H)$-modules is given by:
\beq\lb{tep}
(H\ot_{L(g)}M)\ot (H\ot_{L(h)}N)\cong \oplus_{x \in \cd}(H \ot_{L(\;^xg\;h)}P(x))\\
\eeq
where $\cd$ is a set of representatives for the double cosets $F_g\backslash F/F_h$ and $P(x) \in \mi_{\;^xgh}$ is given by 
\beq
P(x):= (\;^xM\dw^{L(\;^xg)}_{L(\;^xg)\cap L(h)}\ot N\dw^{L(h)}_{L(\;^xg)\cap L(h)})\uw_{L(\;^xg)\cap L(h)}^{L(\;^xgh)}
\eeq
\et 
 The finite group $F$ is defined by $K(H)=kF^*$ and the conjugate modules $\;^xM$ are defined as in Equation \ref{conjm}. \medbreak The above formula shows that the Grothendieck ring structure of $D_{\lam_f}(U,H)$ is also similar to the Grothendieck ring of $D(G)$ of a finite group $G$. The latter ring was described by Cibils in \cite{cibils}. These ring structures were very intensively studied by various authors in connection with the Hochschild cohomology ring and with the rings coming from Green functors. (see for instance \cite{cs}, \cite{bouc}). Witherspoon noticed in \cite{scoh} that abelian cocentral extensions of Hopf algebras have also this type of structure. Connections with the ring construction from \cite{scoh} are given in Section \ref{fus}.
\medbreak
The organization of the paper is the following. In the second section we describe the tools that are needed for the rest of this work.  They include the basic notions of extensions of normal Hopf algebras from \cite{coset} and the results on Clifford theory developed by the author in \cite{coc}.

The third section gives the general description of the irreducible representations of the generalized quantum doubles of $D_{\lam_f}(U,H)$. To this end we prove that $K(U)$ is a normal Hopf subalgebra of $D_f(U,H)$ and apply Clifford's correspondence for this normal Hopf subalgebra.

In Section \ref{cocentral} we prove that a semisimple generalized quantum double $D_{\lam}(U,H)$ is a cocentral extension of a certain Hopf subalgebra by the group algebra of the universal grading group of the category  $\rep(U)$.

Section \ref{catint} discusses the categorical interpretation of the double crossed product. It is shown that the category of representations of a generalized quantum double is the relative center of a certain tensor functor. %Connections with the results from \cite{gnn} are given in this section. 

In the last section we prove the tensor product formula from Theorem \ref{tep}. Using this we show that the Grothendieck ring of a generalized quantum double has a structure similar to the one obtained in \cite{scoh}.
\section{Preliminaries}\lb{prelim}

In this section we recall the basic notions and results on semisimple Hopf algebras and fusion categories that are needed for the rest of the paper.

Throughout of this paper any finite dimensional semisimple Hopf algebra $A$ will be defined over a fixed algebraically closed field $k$ of characteristic zero. Then $A$ is also a cosemisimple Hopf algebra and $S^2=\mtr{Id}$ (see \cite{Lard}), where $S$ is the antipode of $A$. The set of irreducible characters of $A$ is denoted by $\mtr{Irr}(A)$. The Grothendieck group $\mtc{G}(A)$ of the category $\rep(A)$ of finite dimensional left $A$-modules becomes a ring with the multiplicative structure induced by the tensor product of $A$-modules. Then $C(A)=\mtc{G}(A)\ot_{\mathbb{Z}}k$ is a semisimple subalgebra of $A^*$ \cite{Zhc} and it has a basis given by the characters of the irreducible $A$-modules. There is a unique bilinear form  $m_A:C(A)\times C(A)\ra k$ on $C(A)$ determined by the relation $m_A([V],\;[W])=\dim \Hom_A(V,W)$ for any two $A$-modules $V$ and $W$.
%{\bf Definition of the multiplicity}
\subsection{Conjugate modules for normal extensions}
\md
Let $B \subset A$ be a normal Hopf subalgebra of a semisimple Hopf algebra $A$ and let $M$ be an irreducible $B$-module with associated character $\al \in C(B)$. We recall the following notion of conjugate module introduced in \cite{coset}. It was also previously considered in \cite{Schgal} in the cocommutative case.

If $W$ is an $A^*$-module then $W\ot M$ becomes a $B$-module with the following structure:
\bn{equation}\label{def}
 b(w\ot m)=w_0 \ot (S(w_1)bw_2)m
\end{equation}
for all $b \in B$, $w \in W$ and $m \in M$.
Here we used that any left $A^*$-module $W$ is a right
$A$-comodule via $\rho(w)=w_0\ot w_1$.
It can be checked that if $W \cong W'$ as $A^*$-modules then $W\ot M \cong W'\ot M$. Thus for any irreducible character $d \in \mtr{Irr}(A^*)$ associated to a simple $A$-comodule $W$ one can define the $B$-module $\;^dM:= W\ot M$. If $\al \in C(B)$ is the character of $M$ then the character $^d\al$ of $^dM$ is given by

\bn{equation}\label{chfom}
^d\al(x)=\al(Sd_1xd_2)
\end{equation}
for all $x \in B$ (see Proposition 5.3 of \cite{coset}).

\subsubsection{Rieffel's equivalence relations for normal extensions}\lb{rieffe}
%Define two equivalent $B$-modules in the following way. 
We say that two $B$-modules $V$ and $W$ are equivalent if there is an irreducible $A$-module $M$ such that both $V$ and $W$ are irreducible constituents of of $M\dw_B^A$. This defines  an equivalence relation $\sim_B$ on the set $\Irr(B)$ of irreducible characters of $B$ (see \cite{coset}). The same notation is used for the corresponding equivalence relation on $\irr(B)$. There is an analogue equivalent relation on \x \irr(A) \x. Two $A$-modules $M$ and $N$ are equivalent if there is an irreducible $B$-module $V$ that is constituent to both $M_B$ and $N_B$. This equivalence relation is usually denoted by $\sim^A  $. 
\br\lb{corresp}
According to \cite{coset} there is a bijective correspondence between the equivalence classes of $\sim^A  $ and $\sim_B$. This correspondence is given as follows. Let $V$ be an irreducible $A$-module and $M$ a constituent of $V_B$. Then the equivalence class of $V$ under $\sim^A  $ corresponds to the equivalence class of $M$ under $\sim_B$. 
\er
\bn{rem}\label{eqconst}
From Proposition 5.12 of \cite{coset} it follows that the equivalence class of the character $\al \in \mtr {Irr}(B)$ is given by all the irreducible constituents of  $^d\al$ as $d$ runs through all irreducible characters of $A^*$.
\end{rem}

\subsubsection{Restriction and induction to normal Hopf subalgebras}
We recall the following results for induction and restriction to normal Hopf subalgebras from \cite{coset} that will be needed letter.
\bp \lb{ir} Let $B$ be a normal Hopf subalgebra of a \ses \; Hopf algebra $A$. Let $V$ be an irreducible $A$-module with character $\ch \in \irr(A)$ and $M$ be an irreducible $B$-module with associated character $\al \in \Irr(B)$. Suppose that $V$ is a constituent of $M_B$. Then
\beq\lb{restr}
\ch\dw^A_B=\ch(1)\frac{\sum_{\beta\sim_B\al}\beta(1)\beta}{\sum_{\beta\sim_B\al}\beta(1)^2}
\eeq
\beq\lb{ind}
\al\uw^A_B=\frac{\al(1)\dim A }{\dim B}\sum_{\mu\sim^A  \ch}\mu(1)\mu
\eeq
\ep
The two equations above combined give the following:
\beq\lb{indrestr}
\al\uw^A_B\dw^A_B=\frac{\dim A}{\dim B}\frac{\al(1)}{(\sum_{\beta \sim_B \al}\beta^2(1))}\sum_{\beta \sim_B \al}\beta(1)\beta
\eeq
for any $\al \in \Irr(B)$.

Suppose that $\cb$ is an equivalence class of $\sim_B$ corresponding by Remark \ref{corresp} to the equivalence class $\ca$ of $\sim ^A$. Since $A$ is a free left $B$-module it follows that
\beq\lb{restrentr}
\sum_{M \in \ca}(\dim M) M\dw^A_B=\frac{\dim A}{\dim B}\sum_{V \in \cb}(\dim V) V
\eeq

\subsubsection{Subsets closed under multiplication and duality}
Recall from \cite{NR'} that a subset $X \subset \mtr{Irr}(A^*)$ is closed under multiplication if for every $\chi, \mu \in X$ in the decomposition of $\chi\mu=\sum_{\gamma \in \mtr{Irr}(A^*)}m_{\gamma}\gamma$ one has $\gamma \in X$ if $m_{\gamma} \neq 0$. A subset $X \subset \mtr{Irr}(A^*)$ is closed under $``\;^*\;"$ if $x^* \in X$ for all $x \in X$.

Any subset of $X \subset \irr(A^*)$ closed under multiplication and $" \;*\;"$ generates a Hopf subalgebra $A(X)$ of $A$ by 
$$
A(X):=\oplus_{x \in \Irr(A^*)}C_x
$$
where $C_x$ is the simple subcoalgebra of $A$ associated to $x$.

\subsubsection{Definition of the stabilizer of an irreducible module}
Let $M$ be an irreducible $B$-module and $\al \in \mtr{Irr}(B)$ be its character. Then the following result 
was proven in \cite{coc}.
\bn{prop}
The set $\{d \in \mtr{Irr}(A^*)\;| \;^d\al=\eps(d)\al\}$ is closed under multiplication and $``\;^*\;"$. Thus it generates a Hopf subalgebra $Z_A(\al)$ of $A$ that contains $B$.
\end{prop}

\bn{example}
Let $A=kG$ and $B=kN$ for a normal subgroup $N$ of $G$. In this situation $Z_A(\al)$ coincides with the group algebra of the classical stabilizer of $\al$ in $G$ that was introduced by Clifford in \cite{C1}.
\end{example}

In analogy with the group algebra case the Hopf subalgebra $Z_A(\al)$ is called in \cite{coc} the stabilizer of $\al$ in $A$.

\bn{rem}\label{tpr} Let $B$ a normal Hopf subalgebra of a semisimple Hopf algebra $A$ and $M$ be an irreducible $B$-module with character $\al \in \irr(B)$. If $C$ is any subcoalgebra of $A$ then $C\ot M$ has a structure of $B$-module as in Equation \ref{def} using the fact that $C$ is a right $A$-comodule via $\D$. Then it follows that $C\ot M \cong M^{|C|}$ as $B$-modules if and only if $C \subset Z_A(\al)$. \end{rem}
\subsubsection{Clifford theory for Hopf algebras}
Let $B \subset A$ be an extension of normal semisimple Hopf algebras. Fix an irreducible $B$-module $M$ with character $\al \in \Irr(B)$. Suppose that $\mtc{B}_1$ is the equivalence class of $\al$ under Rieffel's equivalence relation $\sim_B$ on $\mtr{Irr}(B)$ for the inclusion $B \subset A$  (see subsection \ref{rieffe}).

\md
As in Remark \ref{corresp} suppose also that $\mtc{A}_1$ is the corresponding equivalence class of $\mtc{B}_1$ on $\mtr{Irr}(A)$. Consider the characters:
\beq
b_1=\sum_{\beta \in \mtc{B}_1}\beta(1)\beta.
\eeq
and
$$a_1=\sum_{\ch \in \mtc{A}_1}\ch(1)\ch.$$

\subsubsection{Rieffel's equivalence relations for the inclusion $B \subset Z_A(\al)$.}
Since $B$ is a normal \hsa $Z_A(\al)$ one can define as above two equivalences relations, on $\mtr{Irr}(Z_{A}(\al))$ respectively $\mtr{Irr}(B)$. %\\Let ${\mtc{Z}_1}, \cdots ,\mtc{Z}_r$ be the equivalence classes in $\mtr{Irr}(Z_{A}(\al))$ and  $\mtc{B'}_1, \cdots ,\mtc{B'}_r$ be the corresponding equivalence classes in $\mtr{Irr}(B)$.

From the definition of the stabilizer it follows that $\al$ by itself form an equivalence class of $\mtr{Irr}(B)$, say $\mtc{B'}_1$. Then as in Remark \ref{corresp} one can let ${\mtc{Z}_1}$ to be the corresponding equivalence class of $\mtc{B}_1$ on \x\irr(Z_A(\al))$ Thus by definition it follows that ${\mtc{Z}_1}$ is given by $${\mtc{Z}_1}=\{\psi \in \mtr{Irr}(Z_{A}(\al))|\;\psi\dw^{Z_A(\al)}_{ _B}
\;\textnormal{contains}\; \al\}.$$

In this situation, for all irreducible characters $\psi \in {\mtc{Z}_1}$ one has by Formula \ref{restr} that $\psi\dw_{ _B}^{Z_{A}(\al)}=\frac{\psi(1)}{\al(1)}\al$.

\subsubsection{Definition of the Clifford correspondence for Hopf algebras}\lb{cc}
As proven in \cite{coc} the above arguments imply that for any $\psi \in \mtc{Z}_1$  all the irreducible constituents of  $\psi\uw^A_{Z_A(\al)}$ are in $\mtc{A}_{1}$. 

\bn{defn} Let $B \subset A$ be a normal extension of semisimple Hopf algebras. We say that Clifford correspondence holds for the extension $B \subset A$  and the irreducible character $\al \in \irr(B)$ if $\psi\uw^A_{Z_A(\al)}$ is irreducible for any irreducible character $\psi \in \mtc{Z}_1$ and the induction function $$\mtr{ind}: {\mtc{Z}_1}\ra \mtc{A}_{1}$$ given by $\mtr{ind}(\psi)=\psi\uw^A_{Z_A(\al)}$ is a bijection.
\end{defn}
\subsubsection{Necessary and sufficient conditions for Clifford correspondence to hold for the character $\al \in \irr(B)$ and inclusion $B \subset A$.}
\md
\bn{thm}\label{main} ( \cite{coc}.) Let $B\subset A$ be a normal extension of semisimple Hopf algebras and $M$ be an irreducible $B$-module with associated character $\al \in \irr(B)$. Then Clifford correspondence holds for $\al$ if and only if  the following equivalent conditions are satisfied:

\bn{enumerate}
\item  $Z_A(\al)$ is a stabilizer in the sense given in \cite{Ri}.
\item The following equality holds:
\beq\lb{must} 
m_{B}(\al, \;\al\uw^A_B\dw^A_B)=m_{B}(\al, \;\al\uw^{Z_A(\al)}_B\dw^{Z_{A}(\al)}_B).
\eeq
\end{enumerate}
\end{thm}

As noted in \cite{Sc1} for a finite group algebra extension $kH \subset kG$ with $H$ normal subgroup of $G$ the above condition is satisfied for any $\al \in \irr(H)$. In this case one obtains the classical Clifford correspondence  for groups from \cite{C1}.
\subsection{Fusion categories associated to semisimple Hopf algebras} \md In this subsection we recall few basics things on fusion categories and then consider the fusion category $\Rep(A)$ associated to a semisimple Hopf algebra $A$. 
\subsubsection{Gradings of fusion categories}\lb{grads}
Let $G$ be a finite group. Recall that a fusion category $\cc$ is $G$-graded if there is a decomposition
\beq\lb{grds}
\cc =\oplus_{g \in G}\cc_g
\eeq
of $\cc$ into a direct sum of full abelian subcategories such that the tensor product of $\cc$ maps $\cc_g \times \cc_h$ to $\cc_{gh}$, for all $g, h \in G$.
\md
There is a universal grading on $\cc$ by a group $U(\cc)$ called the universal grading group of $\cc$. The universal property of this grading consists of the fact that any other grading on $\cc$ is obtained by a quotient group of the universal group $U(\cc)$. \md The fusion category $\cc$ is also called a $G$-extension of $\cc_1$. Gradings and extensions play an important role in the study and classification of fusion categories (see \cite{NG} and \cite{ENO2}). 
\subsubsection{Universal grading group of $\Rep(A)$}
For a semisimple Hopf algebra $A$ let $K(A)$ be the largest central Hopf subalgebra of $A$. Since $K(A)$ is a commutative Hopf algebra and $k$ algebraically closed it follows that $K(A)=kG^*$ for some finite group $G$. \md If $\mtc{C}:=\mtr{Rep}(A)$ is the fusion category of finite dimensional representations of $A$ then it follows that $U(\cc)=G$ (see Theorem 3.8 of \cite{NG}.) Moreover, if
\begin{equation}\label{univ}
\mtc{C}=\oplus_{g \in G}\mtc{C}_g
\end{equation}
is the universal grading of $\rep(A)$ then an irreducible character $\ch$ of $A$ satisfies $\ch \in \mtc{C}_g$ if and only if $\ch\dw^A_{K(A)}=\ch(1)g$. Following the proof of Theorem 3.8 of \cite{NG} it is easy to deduce that
\beq
p_g=\sum_{M \in \co(\cc_g)}e_M
\eeq
where $\co(\cc_g)$ is the set of isomorphism classes of simple object of the abelian category $\cc_g$. Here $e_M\in H$ is the primitive central idempotent associated to $M$.
 \md
On the other hand, by Frobenius reciprocity, it follows that $g\uw^A_{K(A)}=\sum_{\ch \in \cc_g}\ch(1)\ch$. Thus one has  \beq\lb{fpdim}
\sum_{\ch \in \cc_g}\ch(1)^2=\frac{\dim A}{|G|}
\eeq

and therefore $\fp(\cc_g)=\frac{\dim A}{|G|}$ where $\fp(\cc_g)$ is the dimension Frobenius - Peron of $\cc_g$ (see \cite{ENO} for its definition).
%%%%%%%%%%%%%%%%%%%%%%%%%%%%%%%%%
\subsubsection{Nilpotent fusion categories}\lb{nilp}
Recall that a fusion category $\cc$ is called nilpotent if there is a sequence of fusion categories
$$\cc_0 = \mtr{Vec}_k, \cc_1, . . . , \cc_n = \cc$$
such that $\cc_i$ is a $G_i$-extension of $\cc_{i-1}$, for some finite groups $G_1, \cdots ,G_n$.
This is equivalent to the fact that the sequence of adjoint subcategories
\beq
\cc \supseteq \cc_1\supseteq \cc_2  \cdots \cc_n \supseteq \cc_{n+1} \cdots
\eeq
stops at $\cc_m=\mtr{Vec}$ for some $m\geq 1$. Recall from \cite{NG} that the adjoint subcategories are defined  inductively by $\cc_1=\cc_{\mtr{ad}}$ and  $\cc_{n+1}=(\cc_n)_{\mtr{ad}}$ for all $n \geq 1$.
\md
If $\cc=\rep(A)$ for some semisimple Hopf algebra $A$ it follows that $\cc_{\mtr{ad}}=\rep(A//K(A))$. Moreover there are Hopf subalgebras 
\beq
K_1(A)=K(A)\subseteq \cdots K_n(A)\subseteq K_{n+1}(A)\cdots
\eeq
of $A$ such that $\cc_n=\rep(A//K_n(A))$. Thus $\rep(A)$ is a nilpotent fusion category if and only if $K_m(A)=A$ for some $m\geq 1$.
%It follows that $\rep(A)$ is nilpotent if the 
\section{Irreducible modules over generalized quantum doubles}\lb{dct}
In this section a description of the irreducible modules of semisimple  generalized quantum doubles will be given by using Clifford's correspondence for Hopf algebras that was developed in the previous section. 
\subsection{Generalized quantum doubles of Hopf algebras}  Majid introduced in (\cite[Example 7.2.6]{Maj}) the following generalization of  the construction of a Drinfeld double. 
%\subsection{Twisted Hopf algebras}\lb{twsg}
%Let $A$ be a Hopf algebra and $\sigma:A \otimes A \ra k^*$ be an invertible 2-cocycle, that is:
%$$\sigma(x_1, \;y_1)\sigma(x_2y_2,\;z)=\sigma(y_1,\; z_1)\sigma(x, \;y_2z_2)$$
%$$\sigma(1,\;1)=1.$$for all $x, y, z\in A$.
%Define a new Hopf algebra $A_{\sigma}$ with the same comultiplication as $A$ but with the multiplication twisted by %$\sigma$. Thus $$x._{\sg}y=\sg(x_1, y_1) x_2y_2\sg^{-1}(x_3, y_3) $$ for all $x, y \in A$.
%Write the formula for the antipode
Let $U$ and $H$ be Hopf algebras and $\lam:U\ot H\ra k$ be a skew pairing of $U\ot H$. Recall from \cite{Maj} that  this means that $\lam$ is an invertible bilinear map such that
\beq\lb{1}\lam(u, ab)=\lam(u_2, a)\lam(u_1, b)\eeq
\beq\lb{2}\lam(uv, a)=\lam(u, a_1)\lam(v, a_2)\eeq
\beq\lb{3}\lam(1, \;a)=\eps(a)\eeq
\beq\lb{4}\lam(u, \;1)=\eps(u)\eeq
for all $u, v\in U$ and $a, b\in H$.
%check if K(A) is normal in bicrossed product
Let $A:=U\ot H$ be the tensor product Hopf algebra of $U$ and $H$ and let $\sg: A \ot A \ra k$ be the bilinear map defined by 
\beq\lb{cocyc}
\sg(u\ot a,\;v \ot b)=\eps(a)\lam(v, a)\eps(b)
\eeq
for all $u, v\in U$ and $a, b\in H$. Then it can be easily checked that $
\sg$ is an invertible  2-cocycle on $A$ (see also \cite{andf}) with the inverse given $$\sg^{-1}(u\ot a,\;v \ot b)=\eps(a)\lam^{-1}(v, a)\eps(b),$$ for all $u, v\in U$ and $a, b\in H$.
Here $\lam^{-1}(v,\;a)=\lam(Sv, a)=\lam(v, S^{-1}a)$ is the convolution inverse of $\lam$. 
\md
The twist Hopf algebra $A_{\sg}$ is called a crossed product Hopf algebra of $U$ and $H$ with respect to $\lam$ and it is denoted by $D_{\lam}(U,H)$. Note that $U$ and $H$ are both Hopf subalgebras of $D_{\lam}(U,H)$ and the multiplication map $m: U\ot H \ra D_{\lam}(U,H)$ is bijective. Therefore $D_{\lam}(U,H)$ is a particular case of bicrossed product (of $U$ and $H$) in the sense given by Majid in \cite{Maj}.
\md
The multiplication formula in $D_{\lam}(U,H)$ is given by
\beqn
(u\btw a)(v \btw b)=\lam(v_1, a_1)\lam(v_3, S^{-1}a_3)uv_2\bwt a_2b
\eeqn
and the antipode of $D_{\lam}(U,H)$ is given by 
\beqn
S(u\btw a)=S(a)S(u)=\lam(Su_3, Sa_3)\lam(Su_1, a_1)Su_2\btw Sa_2.
\eeqn
\md
These double crossed products $D_{\lam}(U,H)$ are also called generalized quantum doubles since they generalize the well known construction of a quantum double.  %They are also particular cases of bicrossed products of two Hopf algebras.

\subsubsection{Definition of $D_{\lam_f}(U,H)$}
Consider a Hopf algebra map $f : U \ra H^{*\;\mtr{cop}}$ and define $$\lam_f(u, \;h)=f(u)(h)$$ for all $u \in U$ and $h \in H$. Then $\lam_f$ becomes a Hopf skew pairing as above. Moreover the inverse of $\lam_f$ is given by
$$\lam^{-1}_f(u, \;h)=f(u)(S^{-1}h).$$ 
\md
We will also use the notation $<u,h>:=\lam_f(u,h)$. Then the cocycle $\sg$  from Equation \ref{cocyc} becomes 
\beq
\sg(u\ot h, \;u'\ot h')=\eps(u)<u', h>\eps(h').
\eeq
 and the inverse of $\sg$ is given by 
 $$
 \sg^{-1}(u\ot h, \;u'\ot h')=\eps(u)<S^{-1}u', h>\eps(h').
 $$

The above double crossed product Hopf algebra $D_{\lam_f}(U,H)$ will be shortly denoted below by $D_f(U,H)$. 
\bn{rem}
\lb{t}
If $U$ and $H$ are finite dimensional then any skew pairing $\lam$ can be written as $\lam=\lam_f$ for some Hopf algebra map $f : U \ra H^{*\;\mtr{cop}}$.  Indeed, it can be seen that conditions \ref{1}-\ref{4} are equivalent to the map $f:U \ra H^{*\;\;cop}$ given by $u\mapsto \lam(u,\;-)$ to be a bialgebra morphism and therefore a Hopf algebra morphism since both $U$ and $H$ are finite dimensional. Then clearly $\lam=\lam_f$%cite
\end{rem}
\bt\label{central}  Let $f:U\ra H^{*\;\cop}$ be a surjective morphism of semisimple Hopf algebras and $A:=D_f(U,H)$ be the double crossed product as above. Then the Hopf centre $K(U)$ is a normal Hopf subalgebra of $D_f(U,H)$.% {\bf The core of $U$ inside $A$ is $K(U)$.}
\et
%here dont use semismplicity Hopf center exists by Andrusk canad jourrnal
\bpf
Clearly $K(U)$ is closed under the adjoint action of $U$ since it is a central Hopf subalgebra of $U$. It remains to show that $K(U)$ is closed under the adjoint action of $H$ on $A$. 

\md
Identify $A^*$ with $U^* \ot H^{*}$ via 
\beq\lb{idf}
<u^* \ot h^*,\;u\btw h>=<u^*,\;u><h^*,\;h>.
\eeq

If $h \in H$  and $u \in K(U)$ then $$h_1uS(h_2)=u_2 \bwt <u_1,\; h_1><u_3,\;S^{-1}h_3>h_2S(h_4).$$ Then for any $h^* \in H^*$ one has:
\begin{eqnarray*}
(\mtr{id}\ot h^*)  (h_1uS(h_2) ) &=&(\mtr{id}\ot h^*) ( u_2\btw\;<u_1,\; h_1><u_3,\;S^{-1}h_3>h_2S(h_4)) \\
   &=& u_2\btw\; <u_1,\; h_1><u_3,\;S^{-1}h_3><h^*,\;h_2Sh_4>\\
   &=& u_2\btw \;<u_1,\; h_1><u_3,\;S^{-1}h_3><h^*_1,\;h_2><h^*_2,\;Sh_4>\\
  &=& u_2\btw\;<f(u_1)h^*_1(S^{-1}f(u_3))S(h^*_2) ,\;h>\\
  \end{eqnarray*}
Since $u \in K(U)$ and $f$ is a surjective Hopf map it follows that $f(u)\in K(H^{*\;\cop})$. Thus for all $h^* \in H^*$ one has:
\begin{eqnarray*}
(\mtr{id}\ot h^*)  (h_1uS(h_2) ) &=& u_2\btw\; <f(u_1S^{-1}(u_3))h^*_1S(h^*_2) ,\;h>\\
  &=& h^*(1) u_2\btw \;<u_1S^{-1}(u_3) ,\;h>
\end{eqnarray*}
This implies that
\begin{equation}\label{adj1}
   h_1uS(h_2)=<u_1S(u_3) ,\;h>\; u_2\;  \btw1
\end{equation}
for all $h \in H$. Therefore $K(U)$ is also closed under the adjoint action of $H$ on $D_f(U,H)$.
%{\bf $\eps\uw_U^A\dw^A_U$ has the adjoint action?}
\epf

\md
Next Lemma can be regarded as generalization of Theorem 3 from \cite{mintr}.
\bl
Let $f:U\ra H^{*\;\cop}$ be a morphism of semisimple Hopf algebras and $A:=D_f(U,H)$ be the double crossed product as above. Then 
$$D_f(U,H) \cong D_{f^*}(H^{\op}, U^{\op})^{\mtr{op}}$$ as Hopf algebras via $u\ot h \mapsto h \ot u$.
\el
%Radford e la fel
\bpf
Note that $\lam_{f^*}(h,u)=\lam_f(u,h)=f(u)(h)$. Then the rest of the proof is a straightforward computation. 
\epf

\bc
Let $f:U\ra H^{*\;\cop}$ be an injective morphism of semisimple Hopf algebras and $A:=D_f(U,H)$ be the double crossed product as above. Then $K(H)$ is a normal Hopf subalgebra of $A$.
\ec
\bpf

Clearly $f^*:H^{\;\op}\ra U^{*}$ is a surjective morphism of Hopf algebras and one can apply Theorem \ref{central} to $D_{f^*}(H^{\op}, U^{\op})$. Then the previous Lemma proves that $K(H)$ is also a normal Hopf subalgebra of $D_f(U,H)$. Indeed $K(H^{\op})$ is a normal Hopf subalgebra of $D_{f^{*}}(H^{\op},U^{\op})$ and by the above isomorphism this is sent to the Hopf subalgebra $K(H)$ of $D(U,H)$.
\epf
For any finite group $G$ and any $g \in G$ let $C_G(g)$ be the centralizer of the element $g\in G$. We  identify the irreducible $kG^*$-modules $kp_g$ with their characters $g \in G$. 
\subsection{On the central Hopf subalgebra map}
Let $f:U\ra H^{*\;\cop}$ be a surjective morphism of semisimple Hopf algebras. Suppose that $K(U)=kG^*$ and $K(H^{*\;\cop})=kF^*$ for two finite groups $G$ and $F$. Since f is surjective one has that $f(K(U))$ is a Hopf subalgebra of $K(H^{*\;\cop})$. This implies that $$f|_{K(U)}:kG^* \ra kF^*$$ is a morphism of Hopf algebras. Therefore $f|_{K(U)}$  is induced by a group morphism $f^*:F \ra G$.  Then for any $g \in G$ it follows that 
\beq\lb{grm}
f(p_g)=\sum_{\{x \in F\;|\;f^*(x)=g\}}q_x
\eeq
where $\{p_g\}_{ _{g \in G}}$ is the dual group element basis in $kG^*$. Similarly $\{q_x\}_{ _{x \in F}}$ is the dual group element basis in $kF^*$.
 Now let $\cc_1=\rep(U)$ and $\cc_2=\rep(H^{*\;\cop})$ be the fusion categories of finite dimensional representations of $U$ and respectively $H^{*\;\cop}$. Thus as in Subsection \ref{grads} we have universal gradings:
 $$\cc_1=\oplus_{g \in G}(\cc_1)_g$$ and  $$\cc_2=\oplus_{x \in F}(\cc_2)_x.$$ 
 \md
Note that if $K$ is a Hopf subalgebra of $H$ then $U \btw K$ is a Hopf subalgebra of $D_f(U,H)$.\md
 \bt\label{conj}  Let $f:U\ra H^{*\;\cop}$ be a surjective morphism of semisimple Hopf algebras and $A:=D_f(U,H)$ be the double crossed product associated to $f$. Suppose as above that $K(U)=kG^*$ and $K(H^{*\;\cop})=kF^*$.
\bne
\item For any $g \in G$ the conjugate $kG^*$-modules of $g$  in $A$ (in the sense of \ref{def}) are all the group conjugate elements $f^*(x)gf^*(x)^{-1}$ with $x \in F$. 
\item The
stabilizer of a such $kG^*$-module $g \in \Irr(kG^*)$ is given by 
\beq
Z_{A}(g)=U\bwt L(g)
\eeq 
where $L(g)$ is the Hopf subalgebra of $H$ such that
\beq\lb{defl}
\mtr{Rep}(L(g)^*)=\oplus_{\{x\;|\;f^*(x) \in C_G(g)\}}(\mtc{C}_2)_x.
\eeq
\ene
\et

\bn{proof}
\bne
\item
Since $K(U)$ is a central subalgebra in $U$ by relation \ref{chfom} it follows that every irreducible character of $U$ stabilizes the character $g\in G$. 

Let now $d \in \Irr(A^*)$ be an irreducible character of $A^*$. Then applying Equation \ref{adj1} for $S(d)$ one has the following formula for the conjugate character $^dg$:
\begin{eqnarray*}
(^{d}g)(x) : & = & g(Sd_1xd_2)\\ & = & g((Sd)_2xS((Sd)_1)) \\ & = &
<Sx_1x_3,\;S(d)>g(x_2) \\ & = & <Sx_3x_1,\;d>g(x_2).
\end{eqnarray*}
for all $x \in K(U)$. If $x=p_h$ with $h \in G$ then one has that $$(^{d}g)(p_h)=\sum_{uvw=h}<p_{w^{-1}}p_u,\;d>g(p_v)=<\sum_{\{u\;|\;ugu^{-1}=h\}}p_u,\;d>$$
Therefore if $d \in (\cc_2)_x$ then $^dg=\eps(d)f^*(x)gf^*(x)^{-1}$.
\item
From the previous formula it follows that $d \in (\cc_2)_x$ stabilizes $g$ if and only if $f^*(x)g=gf^*(x)$.
%$$\eps(d)g(p_h)= <\sum_{\{u\;|\;ugu^{-1}=h\}}p_u,\;d>$$ for all $h \in G$ or $\delta_{h,g}\eps(d)=<\sum_{\{u\;|\;ugu^{-1}=h\}}p_u,\;d>$.
Therefore the stabilizer $Z_A(g)$ of $g$ is generated by the Hopf subalgebra  $U$ and $L(g)$ where $L(g)$ is the Hopf subalgebra of $H$ with the property:
\beq
\mtr{Rep}(L(g)^*)=\oplus_{\{x\;|\;f^*(x) \in C_G(g)\}}(\mtc{C}_2)_x.
\eeq
\ene
\end{proof}
\bn{rem}\lb{ide}
From the proof above one can see that the irreducible characters  $d \in \irr(L(g)^*)$ are determine by the relation
$$<Sx_3x_1,d>g(x_2)=g(x)\eps(d)$$ for all $x \in K(U)$.
\end{rem}

%\br\lb{com}
%Remark that if $G$ is abelian group then $C_G(g)=G$ and consequently $L(g)=H$. Thus any $g \in G$ is a $D_f(U,H)$-stable character in this case.
%\er
\subsection{Clifford Theory for generalized quantum doubles}
\bn{thm}
Let $f:U\ra H^{*\;\cop}$ be a surjective morphism of semisimple Hopf algebras and $A:=D_f(U,H)$ be the generalized quantum double associated to $f$. Then Clifford's correspondence holds for the normal extension $K(U) \subset D_f(U,H)$ and for any $g \in \Irr(K(U))$.
\end{thm}

\bn{proof}
Let as above $Z_A(g):= U\bowtie L(g)$ be the stabilizer of $g$.
By Theorem \ref{main} to show that Clifford correspondence holds for the extension
$K(U) \subset D_f(U,H)$ one has to prove that
$$m_{K(U)}(g\uw^{Z_A(g)}_{K(U)}\dw^{Z_A(g)}_{K(U)},
\;g)=m_{K(U)}(g\uw^{A}_{K(U)}\dw^{A}_{K(U)}, \;g)$$ for all $g \in G$.

Note that $F$ acts on $G$ via conjugation by $f^*$:
\beq\lb{act}
\;^xg:=f^*(x)gf^*(x)^{-1}.
\eeq Thus the size of the orbit $\mtc{O}(g)$ of $g$ is $\frac{|F|}{|\stab_F(x)|}$ where $\stab_F(x):=\{x \in F\;|\; f^*(x)\in C_G(g)\}$. %Let also be the orbit of $g$ under this action.
 
Theorem \ref{conj} implies that the conjugate $A$-modules of $g$ are exactly the elements of $\co(g)$,  the orbit of $g \in G$ under the action of $F$. Then by Proposition \ref{ir} one has 
\beq\lb{ma}
g\uw^{A}_{K(U)}\dw^{A}_{K(U)}=\frac{\dim A }{\dim K(U)|\co(g)|}\sum_{h \in \mtc{O}(g)}h
\eeq
Since $|\co(g)|=\frac{|F|}{|\stab_F(g)|}$ and $\dim K(U)=|G|$ one has that 
\beq\lb{maf}
m_{K(U)}(g,\;g\uw^{A}_{K(U)}\dw^{A}_{K(U)})=\frac{\dim A}{\dim K(U)|\co(g)|}=\frac{(\dim A)|\stab_F(g)|}{|F||G|}
\eeq 
On the other hand since $g$ is $Z_A(g)$-stable one has by the same Proposition \ref{ir} that $$g\uw^{Z_A(g)}_{K(U)}\dw^{Z_A(g)}_{K(U)}=\frac{\dim Z_A(g)}{\dim K(U)}g.$$ Then it follows
that
\beq\lb{mz}
m_{K(U)}(g,\;g\uw^{Z_A(g)}_{K(U)}\dw^{Z_A(g)}_{K(U)})=\frac{ \dim Z_A(g)}{\dim K(U)}=\frac{(\dim L(g))(\dim U)}{\dim K(U)}.
\eeq

On the other hand formula \ref{fpdim} and definition \ref{defl} of $L(g)$ imply that
 $$\dim L(g)=\frac{\dim H}{|F|}|\{x \in F\;|\; f^*(x)\in C_G(g)\}|=\frac{\dim H}{|F|}|\stab_F(g)|.$$ Thus $$\dim Z_A(g)=\frac{(\dim H)( \dim U)|\mtr{Stab}_F(g)|}{|F|}=\frac{(\dim A) |\mtr{Stab}_F(g)| }{|F|}$$ and Equation \ref{mz} becomes
 \beq\lb{mz'}
m_{K(U)}(g,\;g\uw^{Z_A(g)}_{K(U)}\dw^{Z_A(g)}_{K(U)})=\frac{ \dim Z_A(g)}{\dim K(U)}=\frac{(\dim A)|\mtr{Stab}_F(g)|}{|F||G|}.
\eeq
Comparing Equations \ref{maf} and \ref{mz'} it follows the Clifford correspondence holds for any $g \in \irr(K(U))$.
\epf

For any $g \in G$ define the set 
\beq\lb{i}
\mtc{I}_g=\{ M \in \Irr(U\bwt L(g))\;| M_{K(U)}=g \dim\;M\}
\eeq

\md
Since Clifford correspondence holds for any $g \in \Irr(kG^*)$ one can now give a proof of  Theorem \ref{moddescr}.\md

{\bf Proof of Theorem \ref{moddescr}:} Applying Clifford's correspondence from Subsection \ref{cc} it follows that any irreducible $D_f(U,H)$-module that seats over $g \in G$ is induced from an irreducible module of $U\bowtie L(g)$. Let $M$ be such a module. Then it is easy to show
\beqn
D_f(U,H)\ot_{U\bowtie L(g)}M\cong H\ot_{L(g)}M
\eeqn via the map $$(u \bowtie h)\ot_{U\bowtie L(g)} m \mapsto <u_1, Sh_1><u_3,h_3>h_2\ot_{L(g)} u_2m.\;\;\;\square$$

\md
Denote by $S_{g, M}$ the simple $D_f(U, H)$-module $H\ot_{L(g)}M$ from above. One also has the following:
\md

\bc\lb{bij} With the hypothesis from the previous Theorem there is a bijection $S_{g,M}\mapsto (g, M)$ between the set of irreducible $D_f(U,H)$-modules and the pairs $(g, M)$ with $M \in \ci_g$ and $g \in \Gamma$ where $\Gamma \subset G$ is a set of representative elements of the orbits of the action \ref{act} of $F$ on $G$. 
\ec
\subsubsection{On the Drinfeld double of a Hopf algebra}\lb{da}
Let $A$ be a semisimple Hopf algebra. It is well known that the Drinfeld double $D(A)$ is a double crossed product of $U=A^{*\;\cop}$ with $H=A$ where the map $f : U^{*\;\cop}\ra H$ is the identity $\id_{A^{*\;\cop}}$. 
Therefore $$D_{\id_{A^{*\;\cop}}}(A^{*\;\cop}, A) \xrightarrow{\sim} \;D(A)$$ as Hopf algebras via the identity map. %Thus $U=A^{*\;\cop}$ and $H=A$ in this situation. 
Then the group morphism $f^*$  from Equation \ref{grm} also becomes the identity morphism from $G$ to itself. Thus $K(A^{*\;\cop})=k[G]^{*\;\cop}$ is a normal Hopf subalgebra of $D(A)$. For any element $g \in G_{}$ one has that the stabilizer $Z_{D(A)}(g)$ is given by $A^{*\;\cop}\btw L(g)$. Thus Proposition \ref{moddescr} gives that any irreducible $D(A)$-module is of the form $A^{*\;\cop}\ot_{L(g)}M$ for some irreducible left $A^{*\;cop}\btw L(g)$-module $M \in \ci_g$.
\md
In the example bellow we will see that this description for irreducible modules generalizes the well known description for irreducible modules over $D(G)$, the Drinfeld double of a finite group $G$ (see \cite{L87} or \cite{Sc1}).
%\bn{example}
%Suppose that $H=kG$ and $f:U\ra H^{*\;\cop}$ is a surjective morphism of Hopf algebras. Thus $f:K(U)\ra 
%\end{example}
\subsubsection{The Drinfeld double of a group}
Let $G$ be a finite group and $D(G)$ be the quantum double associated to $G$.
It is well known that any irreducible $D(G)$-module is of the form $kG\ot_{kC_G(g)} M$ for an irreducible $C_G(g)$-module (see for example \cite{W96}.) Next example shows that this description follows from Theorem\ref{moddescr}.

\bn{example}\lb{dg}
Let $G$ be a finite group. Then as above $D(G)\xrightarrow{\sim} \;D_{\id}(kG^{*\;\cop},kG)$
as Hopf algebras. Thus in this situation $K(kG^{*\;\cop})=kG^{*\;\cop}$. Moreover since $\Irr(kG^*)=G$ Formula \ref{defl} implies that the Hopf subalgebra $L(g)$ of $kG$ coincides to the group algebra $kC_G(g)$ of the stabilizer subgroup of $g$. Then it is also easy to check that $\ci_g=\Irr(C_G(g))$ in this situation. Indeed, an irreducible $kG^{*\;\cop}\btw kC_G(g)$ module with $M|_{kG^{*\;\cop}}=g (\dim \;M)$ is completely determined by the action of $C_G(g)$. Thus in this case one obtains the well known description of irreducible modules over $D(G)$ from \cite{W96}. \end{example}
%%%%%%%%%%%%%%%%%
\section{Generalized quantum doubles as cocentral extensions and equivariantizations}\lb{cocentral}
Recall that a group $G$ acts by tensor equivalences on a fusion category $\cc$ if there is a tensor functor $\rho: \underline{G} \to \underline{\mtr{Aut}}_{\otimes} \cC$. This means that, for every $g \in G$, there is a
$k$-linear tensor functor $\rho^g: \cC \to \cC$ and natural isomorphisms of tensor functors
$\rho^{g,h}_2 : \rho^g \rho^h \to \rho^{gh}$, $g, h \in G$, and unity tensor isomorphism
$\rho_0 : \id_\C \to \rho^e$, subject to some natural compatibility  conditions (see \cite{invtamb}). If $G$ acts by tensor equivalences on a fusion category $\cc$ then one can construct the fusion category $\cc^G$ of invariant objects as in \cite{invtamb}.  Simple objects of general equivariantizations will be described in \cite{buna}.
\subsubsection{The equivariant category associated to a cocentral extension}
Recall that an exact sequence of Hopf algebras 
\beq
k \ra A\xrightarrow{i} H\xrightarrow{\pi} kG\ra k
\eeq

is called cocentral if $kG^*\subset \mtc{Z}(H^*)$ via $\pi^*$. Following Proposition 3.5 of \cite{natalecoc} it follows that for any such extension one has that $G$ acts on $\Rep(A)$ and $\Rep(H)=\Rep(A)^G$.
\subsection{The cocentral extension associated to a double crossed product} Let $f :U\ra H^{*\;\cop}$ be any morphism of Hopf algebras and $D_f(U,H)$ be the double crossed product associated to $f$. As in Equation \ref{idf} identify \x D_f(U,H)^* \xrightarrow {\sim} U^* \ot H^*\; \x  as algebras via the evaluation 
\beq\lb{idff}
<u^*\ot h^*, \;u\btw h>=u^*(u)h^*(h)
\eeq
for all $u^* \in U^*$, $u\in U$, $h^*\in H^*$ and $h \in H$.
\bl
Under the identification \x D_f(U,H)^*  \cong U^* \ot H^*\; \x from Equation \ref{idff} the comultiplication on \x D_f(U,H)^* \;\x   becomes the following:
\beq\lb{comult}
\Delta(u^*\ot h^*)=(u_1^*\ot f(e_i)h^*_1f(e_j))\ot (e_i^*u^*_2S^{-1}(e_j^*) \ot h^*_2)
\eeq
where $\{e_i\}_i$ is a basis on $U$ and $\{e_i^*\}_i$ is the dual basis on $U^*$.
\el

\bpf
One has 
\begin{eqnarray*}
(u^*\ot h^*)((x \btw m)(y \btw l))  & = &   (u^*\ot h^*)(xy_2\btw m_2l)<m_1,\; y_1> <m_3, Sy_3>\\ & = & u^*(xy_2)h^*(m_2l)<m_1,\; y_1><m_3, Sy_3> \end{eqnarray*}
 On the other hand: \begin{eqnarray*}
& &<u_1^*\ot f(e_i)h^*_1f(e_j),\;x \btw m><e_i^*u^*_2 S^{-1}(e_j)^* \ot h^*_2, y \btw l> \\ & = & u_1^*(x)(f(e_i)h^*_1f(e_j))(m)(e_i^*u^*_2S^{-1}(e_j^*))(y)h^*_2(l)\\ & = & u_1^*(x)f(e_i)(m_1)h^*_1(m_2)f(e_j)(m_3)(e_i^*)(y_1)u^*_2(y_2)S^{-1}e_j^*(y_3)h^*_2(l)\\ & = &  u_1^*(x)f(y_1)(m_1)f(S^{-1}y_3)(m_3)h^*(m_2l)u^*_2(y_2)\\ & = & u_1^*(xy_2)h^*(l_2m)<m_1,\; y_1><m_3, Sy_3> \end{eqnarray*}
for all $x,y \in U$ and $m,l\in H$.\epf% and $u^*\in U^*$, $h^*\in H^*$.\epf

\bt\lb{cocentral extension} Let $f:U\ra H^{*\;\cop}$ be a Hopf algebra morphism. Suppose that $K(U^*)=kG^*$ and let $\pi_{K(U^*)}:U^*\ra U^*//K(U^*)$ be the canonical Hopf projection. With the above notations it follows that:
\beq
k \ra D_{f\pi_{K(U^*)}^*}((U^*//K(U^*))^*,H)\ra D_f(U,H)\ra kG\ra k
\eeq
is a cocentral extension of Hopf algebras. Therefore $$\rep(D_f(U,H))\cong \rep(D_{f\pi_{K(U^*)}^*}((U^*//K(U^*))^*,H))^G$$ as fusion categories.
\et
\bpf
It is a straightforward computation to verify that the Hopf subalgebra $D_{f\pi_{K(U^*)}^*}((U^*//K(U^*))^*,H)$ is normal inside \x D_f(U,H)$. Moreover $$D_f(U,H)//D_{f\pi_{K(U^*)}^*}((U^*//K(U^*))^*,H) \cong U//((U^*//K(U^*))^*$$ and $$U//(U^*//K(U^*))^*\cong K(U^*)^*=kG.$$

It remains to show that the above sequence is cocentral namely that $K(U^*)$ is a central Hopf subalgebra of $D_f(U,H)^*$. The identification from Equation \ref{idff} shows that $K(U^*)$ is a central subalgebra of $D_f(U,H)^*$ since it is a central subalgebra of $U^*$. It remains to show that $K(U^*)$ is a subcoalgebra of $D_f(U,H)^*$. In order to see this one has to use the comultiplication Formula \ref{comult} from above. Indeed if $u^* \in K(U^*)$ then
\beqarn
\Delta ( u^*\ot \eps_{H^*} ) & = & \sum_{i,j} (u^*_1\ot f(e_i)f(e_j))\ot( e_i^*u_2^*S^{-1}(e_j^*)\ot \eps_{H^*})
\\ &=& \sum_{i,j} (u_1^*\ot f(e_ie_j) )\ot( e_i^*S^{-1}(e_j^*)u_2^*\ot \eps_{H^*}) \\ &=& (u^*_1\ot \eps_{H^*} )\ot(u^*_2\ot \eps_{H^*})
\eeqarn
Here we used that 
\beq
e_i^*S^{-1}(e_j^*) \ot e_ie_j=\eps_{U^*}\ot 1
\eeq
which can be proved by straightforward computation.
\epf

Using the reconstruction theorem from \cite{AD} it follows that 
\beq\lb{andc}
D_f(U,H) \xrightarrow{\sim}D_{f\pi_{K(U^*)}^*}((U^*//K(U^*))^*,H)\;^{\tau}\#_{\sg} \;kG
\eeq
for some cocycle $\sg$ and dual cocyle $\tau$. In the case of a nilpotent fusion category $\rep(U^*)$  (in the sense of \cite{NG}) one can iterate this construction and get the following:
%Recall the following definition of nilpotent fusion categories from \cite{NG}.
\bc
Suppose that $\rep(U^*)$ is a nilpotent fusion category. Then $D_f(U, H)$ is isomorphic as Hopf algebras with the an iteration of cocentral extensions:
$$H \;^{\tau_{r}}\#_{\sg_r}kG_r\;^{\tau_{r-1}}\#_{\sg_{r-1}}\cdots \;^{\tau_{1}} \#_{\sg_1}kG_1.$$
\ec
\bpf
Since $\rep (U^*)$ is a nilpotent fusion category, with the notations from Subsection \ref{nilp} it follows that $K_n(U^*)=U^*$ for some $n \geq 1$. Then one has to iterate for $n$-times the construction from Equation \ref{andc} to get the above isomorphism.
\epf

\section{Categorical interpretation of the double crossed product}\lb{catint}
It is well known (see \cite{Kas}) that $D(A)$-mod is equivalent as braided categories to the center $\cZ(A-\mtr{mod})$ of the category of $A$-modules. We generalize below this example using the relative center of a monoidal functor.

\subsection{The relative center of a monoidal functor}
Let $\cd$ and $\cc$  be monoidal categories, and let $(F, F_2, F_0) : \cd \ra \cc$ be a monoidal functor. Recall that in this situation $F_2(X,Y):F(X\ot Y)\ra F(X)\ot F(Y)$ and $F_0:F(1_{\cd})\ra 1_{\cc}$ are isomorphism satisfying the compatibility axioms from \cite{eno}.

As in \cite{BN2} define the relative center $\cZ_F (\cc)$ of $F$ as the following category: 
\md
Objects of $\cZ_F (\cc)$ are pairs $(V, \gm)$, where $V$ is an object of $\cc$ and $\gm$ is a natural transformation $$\{\gm_X :V\ot F(X) \ra F(X) \ot V\}_{X\in \cd}$$ satisfying:
\beq\lb{hexagon}
( \id_V\ot F_2(X,Y))\gm_{X\ot Y} =(\gm_{X} \ot\id_{F(Y)})(\id_{F(X)}\ot\gm_{Y})(F_2(X,Y)\ot \id_V), 
 \eeq
 and
\beq
 (F_0 \ot \id_V) = \gm_{1_{\cd}}(\id_V \ot F_0)
 \eeq
 \md
 Morphisms $u:(V, \gm)\ra (W, \delta)$  in $\cZ_F (\cc)$ are morphisms $u:V\ra W$ in $\cc$ that are compatible with $\gm$ and $\delta$, i.e 
 \beq
 (u \ot \id_X)\gm_X=\delta_X(\id_X \ot u)
 \eeq
for any $X \in \cd$.
The relative center category $Z_F (\cc)$ is a monoidal category (see \cite{BN2}), with the tensor product defined by
 \beq
(V, \gm) \ot (W, \gm') = (V \ot W, (\gm \ot \id_W)(\id_V \ot \gm')).
 \eeq
 Moreover there is a canonical forgetful functor $U : \cZ_F (\cc) \ra \cc$ which is clearly monoidal.
 \br
Note that $\cc$ becomes a $\cd$-bimodule category via $F$ and the relative center $Z_F(\cc)$ from above coincides to the relative center of a bimodule category  introduced in \cite{gnn}.
\er

There are also two canonical functors 
\beq\lb{cic}
F_{\cd} : \cZ(\cd) \ra \cZ_F (\cc)
\eeq
and 
\beq\lb{cid}
F_{\cc}:\cZ(\cc)\ra \cZ_F (\cc)
\eeq obtained from the braidings of $\cZ(\cd)$ and $\cZ(\cc)$ as following.
The functor $F_{\cc}$ is obtained by restriction of the braiding on $\cZ(\cc)$ to the image of $F$. The other functor $F_{\cd} : \cZ(\cd) \ra \cZ_F (\cc)$ is defined by applying $F$ to the braiding from $\cZ(\cd)$.
\br
Note that for $\cc=\cd$ and $F$ identity one obtains the usual Drinfeld center $\cZ(\cd)$ of the fusion category $\cd$.
\er

The category $\cZ_{F}(\cc)$ is also denoted by $\cZ_{\cd}(\cc)$ if the functor $F$ is implicitly understood.
\subsection{Relative centers and double crossed products}
Suppose that $f:U\ra H$ is a morphism of semisimple Hopf algebras. Therefore $f$ induces a functor $f_{*}:\rep(H) \ra \rep(U)$ and one can consider the relative center $\cZ_{f_{*}}(\rep(H))$.

The proof of the following Proposition follows the lines from \cite{Kas}.
\bp\lb{interp}
With the above notations one has the following equivalence of fusion categories:
\beq
\rep(D_{f^{*\;\cop}}(H^{*\;\cop},U))\cong \cZ_{f_{*}}(\rep(H))
\eeq
\ep

\bpf
Let $\cc:=\rep(U)$ and $\cd:=\rep(H)$. Thus one has $f_{*}:\cd \ra \cc$. If $(V, \gm^V)\in \cZ_{f_{*}}(\cc)$ then define $\rho(v):=\gm^V_{H}(1\ot v)\in V\ot H$.  The property $$\gm^V_{H\ot H}=(\gm^V_{H} \ot \id)(\id \ot \gm^V_{H})$$
implies that $V$ is a $H$-right comodule and therefore a left  $(H)^{*\;\cop}$ - module. On the other hand the fact that $\gm^V_{H}$ is $H$-linear implies that the left $H$-action and the left $U$-action satisfy the necessary compatibility conditions in order to obtain that $V$ is a $D_{f^{*\;\cop}}(H^{*\;\cop},U)$-module. In this way one obtains a tensor functor $F: \cZ_{f_{*}}(\cc)\ra \rep(D_{f^{*\;\cop}}(H^{*\;\cop},U))$ given by $(V, \gm^V)\mapsto V$ with the above $D_{f^{*\;\cop}}(H^{*\;\cop},U) $-module structure of $V$.

Conversely, define the functor $G: \rep(D_{f^{*\;\cop}}(H^{*\;\cop},U))\ra\cZ_{f_{*}}(\cc)$ as following. For any $D_{f^{*\;\cop}}(H^{*\;\cop},U)$-module $V$ define $\gm^V_X:X\ot V\ra V\ot X$ by $x \ot v \mapsto v_0\ot v_1x$ where  $\rho(v)=v_0\ot v_1$ is the $H$-comdule structure induced on $V$ by the left $H^{*\;\cop}$-module structure.  It can be checked that $\gm^V_X$ is $U$-linear and the diagram  \ref{hexagon} is satisfied. Moreover $\gm^V_X$ is invertible with the inverse given by $x \ot v \mapsto v_0\ot S(v_1)x$.
\epf

Note that there is a Hopf algebra map $D_{f^{*\;\cop}}(H^{*\;\cop},U)\ra D(U)$ given by $h^* \btw u\mapsto f^{*\;\;\\cop}(h^*)\btw u$. At the category level this map induces the functor  ${f_{*}}_{\cc}:\cZ(\cc)\ra \cZ_{f_{*}}(\cd)$ from \ref{cid}. Also the Hopf algebra map $D_{f^{*\;\cop}}(H^{*\;\cop},U)\ra D(H)$ given by $h^* \btw u\mapsto h^*\btw f(u)$ induces the other functor ${f_{*}}_{\cd}:\cZ(\cd)\ra \cZ_{f_{*}}(\cd)$ from \ref{cic}.
\subsection{On the quantum double of a semisimple Hopf algebra.}% On the results from \cite{gnn}}
Suppose that $A$ is a semisimple Hopf algebra and let $D(A)$ be its quantum double. Then as we already remarked $D(A)$ is a double crossed product of $A$ and $A^{*\;\cop}$, namely $D(A)=D_{\id}(A^{*\;\cop}, A)$. Therefore one can apply the results of Theorem \ref{cocentral extension} section to this particular case. %We will show that in this case we will obtain the results from \cite{gnn}.
\md
Let $G=U(A)$ be the universal grading group of the fusion category $\cc=\rep(A)$. Then the universal grading of $\cc$ is given by
\begin{equation}
\mtc{C}=\oplus_{g \in G}\mtc{C}_g
\end{equation}
where $\cd:=\cc_1=\rep(A//K(A))$. Then as in  Proposition \ref{interp} one obtains that $$\cZ_{\cd}(\cc)=\rep(D_{\pi_{K(A)}^*}((A//K(A))^{*\;\cop}, A))$$ where $D_{\pi_{K(A)}^{*\;\cop}}((A//K(A))^*, A)=(A//K(A))^{*\;\cop}\btw A$.

For the rest of this subsection let $B:= (A//K(A))^{*\;\cop}\btw A$.
\md
\bp The Hopf algebra $K(A)$ is a central Hopf subalgebra of $B$ and therefore there is a canonical induced grading by $G$ on $\rep(B)$. 
\ep

\bpf
Indeed if $x \in K(A)$ and $f \in (A//K(A))^*$ then one has that $xf=(x_1\rh f \lh Sx_3)\ot x_2=f \btw \eps(Sx_3x_1)x_2=f \btw x$. Thus $K(A)$ is central Hopf subalgebra of $B=(A//K(A))^*\btw A$. The induced canonical grading by $G$ is given by $M \in \rep(B)$ if and only if $M\dw^{B}_{K(A)}=(\dim M)g$.
\epf
Note that from the proof of Proposition \ref{interp} it follows that the graded components are given by $\rep(B)_g=\cZ_{\cd}(\cc_g)$.
Thus we obtain the grading 
\beqn
\cZ_{\cd}(\cc)=\oplus_{g \in G}\cZ_{\cd}(\cc_g)
\eeqn
which coincides with the grading $(18)$ from page 10 of \cite{gnn}.
\md
On the other hand applying Theorem \ref{cocentral extension} it follows that one has the following cocentral extension of $D(A)$:
\beq
k \ra B \ra D(A)\ra kG\ra k
\eeq

Using Proposition 3.5 of \cite{natalecoc} this proves the fact $\cZ_{\cd}(\cc)^G\cong \cZ(\cc)$.  Moreover using the notation from Section \ref{dct} it follows that $\ci_g$ coincides to the set of simple $C_G(g)$-equivariants objects of $\cZ_{\cd}(\cc_g)^{C_G(g)}$. Then the bijection from Corollary \ref{bij} corresponds to the bijection described in Proposition 3.9 of \cite{gnn}.

\section{On the Grothendieck ring of a generalized quantum double}\lb{fus}
\subsection{Hopf subalgebras arising from the dual universal grading} Let $H$ be a semisimple Hopf algebra with $K(H^*)=kF^*$ for some finite group $F$. Let $\cc=\rep(H^*)$ be the category of finite dimensional left $H^*$-modules and 
\beq\lb{dgr}
\cc:=\oplus_{x \in F}\cc_x
\eeq
be the universal grading of $\cc$. 
\md

\subsubsection{} To any irreducible character $d \in \Irr(H^*)$ one can associate a simple subcoalgebra $C_d$ of $H$ as in \cite{Lar}. Then for any subset $X \subset F$  one can define a subcoalgebra $H(X)$ of $H$ as the sum of the simple subcoalgebras $C_d$ of $H$ verifying:
\beq
H(X)^*=\bigoplus_{\{d \in \Irr(H^*)\;|\;d \in \cc_x \; \text{for some}\; x \in X\}}C_d
\eeq
%Similarly define the Hopf subalgebra $H(N)$ of $H$ generated by $\Irr(H(M)^*)=\Lam_{\cc(N)}$.
%\ncm{\deg}{\mtr{deg}}

\md
By \cite{NZ} it follows $H(X)$ is a Hopf subalgebra of $H$ if and only if $X$ is a subgroup of $F$. Moreover note that
\beq\lb{dimz}
\dim H(X)=|X|\fp(\cc_1)
\eeq
for any subset $X \subset F$.
\md
From Equation \ref{defl} it follows that $L(g)=H(F_g)$ for any $g \in G$.
\md
Suppose that $M$ and $N$ are subgroups of $F$ and $d \in \Irr(H^*)$ with $d \in \cc_x$ under the grading from Equation \ref{dgr}. Note that the double coset $H(M)C_dH(N)$ in $H(M) \backslash H \slash H(N)$  is given by
\beq
H(M)C_dH(N)=\oplus_{t \in MxN}C_t
\eeq 
where $d \in \cc_f$. See \cite{cos} for the definition of double cosets for Hopf subalgebras. Thus there is a bijection between the Hopf algebra double cosets $H(M) \backslash H \slash H(N)$ and the group double cosets $M\backslash  F/ N$ . The coset $H(M)C_dH(N)$ corresponds to the double coset of $MxN \in M\backslash  F/ N$ where $x \in F$ is chosen by $d \in \cc_x$.
\md
\br\lb{mst}
Note that if $N$ is a subgroup of $F$ and $X \subset F$ such that $NX\subseteq X$ then $H(X)$ is a left $H(N)$-module by Theorem of \cite{NZ}. Similarly, if $XN\subseteq X$ then $H(X)$ is a right-$H(N)$ module.
\er
\subsubsection{}
Since $K(H^*)=kF^* \subset H^*$ note that for any $a \in H(X)$ one has that $p_y(a)=0$ if $y \notin X$. On the other hand since $\sum_{x \in F}p_x=\eps_H$ it follows that
\beq\lb{eps}
\sum_{x \in X}p_x(h)=\eps(h)
\eeq
for all $a \in H(X)$ and any subset $X \subset F$.
\md
Note also that for any $x,y \in F$ and $a \in H$
\beq\lb{mult}
p_x(a_1)p_y(a_3)a_2=\delta_{x, y}p_x(a_1)a_2
\eeq
since $p_g$ and $p_h$ are central elements of $H^*$. Indeed for any $f \in H^*$ one has
 \beqarn
f(p_x(a_1)p_y(a_3)a_2) & = & (p_xfp_y)(a) =  \delta_{x, y}(p_xf)(a)\\ & = & f(\delta_{x, y}p_x(a_1)a_2).\eeqarn
\subsection{On the Grothendieck ring structure of $D_f(U,H)$}
As above consider $f:U \ra H^{\;*\;\cop}$ be a surjective Hopf algebra map and let $A:=D_f(U,H)$. Then $F$ acts on $G$ as in Equation \ref{act} via $\;^xg=f^*(x)gf^*(x)^{-1}$. For any $g \in G$ let $F_g$ be the stabilizer of $G$ under this action and denote by $\co(g)$ the orbit generated by $g$. Thus 
\beq
|\co(g)|=\frac{|F|}{|F_g|}
\eeq
%\textcolor[rgb]{.00, .00, 1.00}{
\bl\lb{indc}
Let $i:L\hookrightarrow H$ be a Hopf subalgebra of $H$ and suppose that $M$ is a $D_{i^{*\;\cop}f}(U,L)$-module. Then $H\ot_LM$ is $D_f(U,H)$-module with the regular action of $H$ and the structure of $U$-module given by
\beq\lb{minduh}
u(h\ot_Lm)=<u_1, Sh_1><u_3, h_3>(h_2\ot_Lu_2m)
\eeq
\el
\bpf
Note that $D_{i^{*\;\cop}f}(U,L)=U\btw L \subset U\btw H$ is a subalgebra of $D_f(U,H)$. Then
\beq
M\uw^{D_f(U,H)}_{D_{i^{*\;\cop}f}(U,L)}=(U\btw H)\ot_{U\btw L}M\cong H\ot_LM
\eeq
and the module structure coincides with the one from above.
\epf
\subsubsection{On the conjugate modules $\;^xM$}
%The action of $G$ on $\cup_{g \in G}\mi_g$.}
Let $\mtr{I}_g$ be the abelian full subcategory of $U\btw L(g)$-mod generated by the set of irreducible modules $M \in \ci_g$. %It follows that $S_{g, M}$ is a one dimensional module if and only if $M$ is a one dimensional $U$-module and $\mtr{Im}(f^*)\subset C_G(g)$. 
We extend the notation $S_{g, M}$ for any $M \in \mi_g$, not necessarily irreducible $U\btw L(g)$-module. Thus if $M \in \mi_g$ then 
\beq
S_{g, M}:=H\ot_{L(g)}M
\eeq
with the action of $U$ defined as in Equation \ref{minduh}.
%Let $\mi_g$ be the abelian subcategory of $U \bwt L(g)$-mod generated by the objects of $\ci_g$.
%Fix a subset $\Gamma\subset G$ of representatives of the orbits of the action of $F$ on $G$. 
%For any $x \in F$ and $M \in \ci_g$ there is a module$\;^xM$ such that $$H\ot_{L(g)}M\cong H\ot_{L(\;^xg)}\;^xM$$ as $D_f(U, H)$-modules.
%Thus there is an action of the group $G$ on the abelian category  $\oplus_{g \in G}\ci_g$ given by $\rho^{x}(M)=\;^xM$.
\md

For any $x \in F$ and $M \in \ci_g$ define the following left $L(\;^xg)$-module
\beq\lb{conjm}
\;^xM:=H(xF_g)\ot_{L(g)}M.
\eeq
Since $F_{\;^xg}=xF_gx^{-1}$ it follows $F_{\;^xg}(xF_g)=xF_g$ and therefore $\;^xM$ is a $L(\;^xg)$-module via left  multiplication. Note also that $\dim \;^xM=\dim M$.

Then it can be easily checked that $\;^xM$ is a $U\btw L(\;^xg)$-module via 
\beq
(u\btw a)(h\ot_{L(g)}m)=<u_1,S(a_1h_1)><u_3, a_3h_3>a_2h_2\ot_{L(g)}u_2m
\eeq  
 Indeed applying Lemma \ref{indc} one has that $H\ot_{L(g)}M$ is a $D_f(U,H)$-module and in particular a $U \btw L(\;^xg)$-module. Then it can be checked that $\;^xM$ is a $U \btw L(\;^xg)$-submodule of $H\ot_{L(g)}M$.
 
\bl\lb{xg}
With the above notations it follows that $\;^xM \in \mi_{\;^xg}$. 
\el
\bpf
One has to check that $\;^xM\dw^{\;U}_{\;K(U)}=\;^xg\;\dim M$. Since $M\dw^{\;U}_{\;K(U)}=g\;\dim M$ one has that $p_bm=\delta_{b,g}m$ for all $b \in G$.

If $l \in G$ using the module structure from Lemma \ref{indc} one has
\beqarn
p_l(h\ot_{L(g)}m) & = & \sum_{abc=l}<p_a,h_1><p_c, Sh_3>h_2\ot_{L(g)}p_bm
\\ &= & \sum_{abc=l}<p_a,h_1><p_c, Sh_3>h_2\ot_{L(g)}\delta_{g, b}m
\\ &= & \sum_{agc=l}<p_a,h_1><p_c, Sh_3>h_2\ot_{L(g)}m
\\ &= & \sum_{agc=l}<p_a,h_1><p_{c^{-1}}, h_3>h_2\ot_{L(g)}m
\eeqarn
Applying now Equation \ref{mult} (for $H=U^*$) one has
\beqarn
p_l(h\ot_{L(g)}m) & = & \sum_{agc=l}<p_a,h_1><p_{c^{-1}}, h_3>h_2\ot_{L(g)}m \\ &= & \sum_{agc=l}\delta_{a, c^{-1}}<p_a, h_1>h_2\ot_{L(g)}m
\\ &= & \sum_{\{a \in G\;|\;aga^{-1}=l\}}<p_a, h_1>h_2\ot_{L(g)}m
\eeqarn
On the other hand by formula \ref{grm} it follows that 
\beqarn
\sum_{\{a \in G\;|\;aga^{-1}=l\}}<p_a, h_1>h_2 & = & \sum_{\{a \in G\;|\;aga^{-1}=l\}}\sum_{\{y \in F\;|\; f^*(y)=a\}}q_y(h)
\\ &= & \sum_{\{y \in F\;|\; f^*(y)gf^*(y)^{-1}=l\}}q_y(h)
\eeqarn

Using also Equation \ref{eps} for $X=xF_g$ it follows that $$\sum_{\{a \in G\;|\;aga^{-1}=l\}}<p_a, h_1>h_2 = \sum_{\{y \in xF_g\;|\; f^*(y)gf^*(y)^{-1}=l\}}q_y(h_1)h_2$$

On the other hand it can be easily checked that the following set $\{y \in xF_g\;|\; f^*(y)gf^*(y)^{-1}=l\}$ is not empty if and only if $l=\;^xg$ and in this case it coincides to $xF_g$. Thus 
\beqarn
p_l(h\ot_{L(g)}m) & = & \delta_{l, \;^xg}(\sum_{y \in xF_g}q_y(h_1)h_2)\ot_{L(g)}m= \delta_{l, \;^xg}\;h\ot_{L(g)}m
\eeqarn
which finishes the proof.\epf
\bl\lb{orbit}
With the above notations one has that $H\ot_{H(F_g)}M \cong H\ot_{H(F_{\;^xg})}\;^xM$ as $D_f(U,H)$-modules.
\el
\bpf
One can check that the map 
\beqn
H\ot_{H(F_{\;^xg})}(H(xF_g)\ot_{H(F_g)}M)\xra{\psi} H\ot_{H(F_g)}M
\eeqn
given by 
\beq
a\ot_{L(\;^xg)}(b \ot_{L(g)}m)\mapsto ab\ot_{L(g)}m
\eeq
is an isomorphism of $D_f(U, H)$-modules. Indeed, since $F_{\;^xg}=xF_gx^{-1}$ it follows by Remark \ref{mst} that $\psi$ is a well defined linear map. On the other hand, clearly $\psi$ is a surjective morphism of $D_f(U,H)$-modules. Since both the domain and codomain have the same dimension it follows that $\phi$ is isomorphism.
\epf

\bl\lb{xgy}
Suppose that $M \in \ci_g$. With the above notations one has that:
\beq
\;^{\;^y}(\;^xM)\cong \;^{yx}M
\eeq
as $U \btw L(\;^{yx}g)$-modules.
\el
\bpf
One can check that the following map:
\beq
\;^y(\;^xM)=H(yF_{\;^xg})\ot_{H(F_{ _{\;^xg})}}(H(xF_g)\ot_{H(F_g)}M)\xra {\phi} H(yxF_g)\ot_{H(F_{ _{\;^{yx}g}})}M
\eeq
given bx $\phi(a\ot_{H(F_{ _{\;^xg})}}(b\ot_{H(F_g)}m))=ab\ot_{H(F_{ _{\;^{yx}g}})}m$ is a well defined surjective morphism of $U\btw L(\;^{yx}g)$-modules. Since both modules have the same dimension it follows that this is an isomorphism. Note that $\phi$ is just the restriction of the above map $\phi$ to the domain $H(yF_{\;^xg})\ot_{H(F_{ _{\;^xg})}}(H(xF_g)\ot_{H(F_g)}M)$.\epf

\md
Since $\;^1M=M$ this implies that $\;^xM$ is also an irreducible $U\btw L(\;^xg)$-modules if $M$ is an irreducible $U \btw L(g)$-module. Thus one has that $\;^x(\ci_g)=\ci_{\;^xg}$ and in this way one can define a function 
\beq\lb{cxg}
c_{x,g}:\ci_g \ra \ci_{\;^xg}
\eeq
This define an action of $F$ on $\sqcup_{g \in G}\ci_g$.\md 
\subsubsection{}Note that $L(g)\cap L(h)\subseteq L(gh)$ for any $g,h \in G$. 
\bl\lb{gh}
There is a well defined map:
\beq
m_{g,h}:\mi_g \times \mi_h \ra \mi_{gh}
\eeq
given by 
\beq\lb{mgh}
m_{g,h}(M, N):=(M\dw^{L(g)}_{L(g) \cap L(h)}\ot N\dw_{L(g) \cap L(h)}^{L(h)})\uw^{L(gh)}_{L(g) \cap L(h)}
\eeq
for any $M \in \mi_g$ and $N \in \mi_h$.
\el
\bpf
Indeed $M \in U\btw L(g)$-mod and $N \in U\btw L(h)$-mod then clearly $m_{g,h}(M,N)\in U \btw L(gh)$-mod. Note that the above $U\btw L(gh)$-module structure is obtained by induction as in Lemma \ref{indc}.
Moreover, similarly to Lemma \ref{xg} it can be checked by a straightforward computation that 
$m_{g,h}(M, N)\in \mi_{gh}$. \epf
%=(M\dw^{L(g)}_{L(g) \cap L(h)}\ot N\dw_{L(g) \cap L(h)}^{L(h)})\uw^{L(gh)}_{L(g) \cap L(h)} 
\subsection{On the fusion rules of generalized quantum doubles}
Using a Mackey type argument we obtain the following formula for the restriction of any irreducible $D_f(U,H)$-module to the Hopf subalgebras $U\btw L(h)$ of $D_f(U,H)$.

\bl\lb{mackey}
For any $M \in \ci_g$ one has that
\beq\lb{mck}
H\ot_{\;L(g)}M \cong \bigoplus_{x \in \cd} (L(h)\ot_{L(h)\cap L(\;^xg)}\;^xM\dw^{L(\;^xg)}_{L(h)\cap L({ \;^xg})})%(H(xF_g)\ot_{\;L(g)}M)
\eeq
as $U \btw L(h)$-modules where $\cd$ is a set of representatives for the double cosets $F_h\backslash F/F_g$
\el
\bpf
Note that for any $x \in \cd$ the map
 \beq\lb{ons}
\phi_x:  H(F_h)\ot_{H(F_h)\cap H(F_{ \;^xg})}(H(xF_g)\ot_{\;H(F_g)}M)\ra H(F_hxF_g)\ot_{H(F_g)}M
 \eeq
given by 
 \beq
 a \ot_{\;H(F_h)\cap H(F_{ \;^xg})}(b\ot_{\;H(F_g)}m)\mapsto ab\ot_{\;H(F_g)}m
 \eeq
 is well defined and a surjective morphism of left $L(h)$-modules.  By a direct computation it can be checked that $\phi_x$ is also a morphism of $U$-modules and therefore a morphism of $U\btw L(h)$-modules. Since 
 \beq
 |F_hxF_g|=\frac{|F_h||F_{\;^xg}|}{|F_h\cap F_{ \;^xg}|}
 \eeq it follows by Equation \ref{dimz} that both $U \btw L(h)$-modules have the same dimension and therefore $\phi_x$ is an isomorphism. Then $\phi_{g, h}:=\oplus_{x \in \cd}\phi_x$ provides an isomorphism between the two left $U \btw  L(h)$-modules from Equation \ref{mck}.
\epf

\bl\lb{prind} Suppose that $K$ and $L$ are Hopf subalgebras of $H$ and $M$ and $N$ are left $K$ and respectively $L$-modules. Then
%\beqN\uw^H_L \ot M\uw^H_K \ra (N\uw_L^H\dw^H_K\ot M)\uw^H_K\eeq
\beq
\phi: (H\ot_LN)\ot (H\ot_KM)\ra H\ot_K((H\ot_LN)\dw^H_K\ot M)
\eeq
given by $(a\ot_Ln)\ot (b \ot_Km)\mapsto b_2\ot_K((S^{-1}b_1a\ot_Ln)\ot m)$ is an isomorphism of $H$-modules.
\el
\bpf
It can be checked by straight forward computations that $\phi$ is an isomorphism of $H$-modules with inverse given by $b\ot_K((a\ot_Lm)\ot n)\mapsto (b_1a\ot_Ln)\ot(b_2\ot_Km)$.
\epf

%\red{
\vskip 0.3 cm
{\bf Proof of Theorem \ref{tp}:} Lemma \ref{gh} shows that $P(x)\in \mi_{\;^xgh}$ and therefore the right hand term of Equation \ref{tep} is a $D_f(U,H)$-module.

%\blue{
It can be easily check that the isomorphism 
\beq
M\uw^{\;H}_{\;L(g)}\ot N\uw^{\;H}_{\;L(h)}  \xra{\phi}  (M\uw_{\;L(g)}^{\;H}\dw_{\;L(h)}^H\ot N) \uw_{\;L(h)}^{\;H}
\eeq
from Lemma \ref{prind} is in our situation an isomorphism of $D_f(U,H)$-modules.
Indeed one has 
\beqarn
& &\phi(  u[(a\ot_Ln)\ot (b \ot_Km)]) =  \phi(u[(a\ot_Ln)\ot (b \ot_Km)])
\\ &= & \phi((a_2\ot_Lu_2n)\ot (b_2 \ot_K(u_4)_2m))
\lam(u_1, Sa_1)\lam(u_3, a_3)\\ & & \lam((u_4)_1, Sb_1)\lam((u_4)_3, b_3)
\\ &= &  \phi((a_2\ot_Lu_2n)\ot (b_2 \ot_Ku_5m))
\lam(u_1, Sa_1)\lam(u_3, a_3)\lam(u_4, Sb_1)\lam(u_6, b_3)
\\ &= & (b_2)_2\ot_K((S^{-1}((b_2)_1a_2)\ot u_5m) \lam(u_1, Sa_1)\lam(u_3, a_3)\lam(u_4, Sb_1)\lam(u_6, b_3)
\\ &= & b_3\ot_K((S^{-1}b_2a_2\ot_Lu_2n)\ot u_5m)
\lam(u_1, Sa_1)\lam(u_3, a_3)\lam(u_4, Sb_1)\lam(u_6, b_4)
\eeqarn
On the other hand,
\beqarn
& & u\phi((a\ot_Ln)\ot (b \ot_Km)) = u(b_2\ot_K((S^{-1}b_1a\ot_Ln)\ot m))
\\ &=& (b_2)_2\ot_Ku_2((S^{-1}b_1a\ot_Ln)\ot m)     \lam(u_1, S(b_2)_1)\lam(u_3, (b_2)_3)
\\ &=& b_3\ot_K[u_2(S^{-1}b_1a\ot_Ln)\ot u_3m]       \lam(u_1, S(b_2))\lam(u_4, b_4)
\\ &=& b_3\ot_K[((S^{-1}b_1a)_2\ot_L(u_2)_2n)\ot u_3m]  \\ & &    \lam(u_1, S(b_2))\lam(u_4, b_4)\lam((u_2)_1, S((S^{-1}b_1a)_1))\lam((u_2)_3, (S^{-1}b_1a)_3))
\\ &=& b_3\ot_K((S^{-1}b_2a_2\ot_Lu_3n)\ot u_5m)  \lam(u_1, Sb_4)\lam(u_6, b_6)\\ & & \lam(u_2, S(S^{-1}b_3a_1))\lam(u_4, S^{-1}b_1a_3)
\\ &=& b_3\ot_K((S^{-1}b_2a_2\ot_Lu_3n)\ot u_5m) \\ & &  \lam(u_1, Sb_4)\lam(u_2, S(a_1)b_3)\lam(u_4, S^{-1}b_1a_3)\lam(u_6, b_6)
\\ &=&  b_3\ot_K((S^{-1}b_2a_2\ot_Lu_2n)\ot u_4m)  \\ & & \lam(u_1, S(a_1)b_3Sb_4)\lam(u_3, S^{-1}b_1a_3)\lam(u_5, b_6)
\\ &= & b_3\ot_K((S^{-1}b_2a_2\ot_Lu_2n)\ot u_5m)\\ & &\lam(u_1, Sa_1)\lam(u_3, a_3)\lam(u_4, Sb_1)\lam(u_6, b_4)
%\\ &=& \\ & &\\ &=& \\ & &\\ &=& \\ & &\\ &=& \\ & &\\ &=& \\ & &\\ &=& \\ & &
\eeqarn
Following then Lemma \ref{mackey} it follows that the left hand side of Equation \ref{tep} is
%\blue{
\bn{eqnarray*}
 & & M\uw^{\;H}_{\;L(g)}\ot N\uw^{\;H}_{\;L(h)}  \cong  (M\uw_{\;L(g)}^{\;H}\dw_{\;L(h)}^H\ot N) \uw_{\;L(h)}^{\;H}\\ &\cong & \oplus_{x \in \cd}((L(h)\ot_{L(h)\cap L( \;^xg)}(\;^xM\ot N))\uw_{\;L(h)}^{\;H})
\\ & \cong & \oplus_{x \in \cd}\;H\ot_{{L(\;^xg)\cap L(h)}}(\;^xM\dw^{L(g)}_{L(\;^xg)\cap L(h)}\ot N\dw^{L(h)}_{L(\;^xg)\cap L(h)})
\end{eqnarray*}
as $D_f(U,H)$-modules where the $D_f(U,H)$-module structure of the last term is given by Lemma \ref{indc}.
%}

%\blue{
On the other hand, the $D_f(U,H)$-module $H\ot_{L(\;^xgh)}P(x) $ can be written as
\bn{eqnarray*}
& & H\ot_{L(\;^xgh)}(L(\;^xgh)\ot_{{{L(\;^xg)\cap L(h)}}}(\;^xM\dw^{L(g)}_{L(\;^xg)\cap L(h)}\ot N\dw^{L(h)}_{L(\;^xg)\cap L(h)}))\\ & \cong & H\ot_{{L(\;^xg)\cap L(h)}}(\;^xM\dw^{L(g)}_{L(\;^xg)\cap L(h)}\ot N\dw^{L(h)}_{L(\;^xg)\cap L(h)}).
\end{eqnarray*}
It is also straightforward to verify that under the last isomorphism the $D_f(U,H)$-module structure of 
\beq
H\ot_{{L(\;^xg)\cap L(h)}}(\;^xM\dw^{L(g)}_{L(\;^xg)\cap L(h)}\ot N\dw^{L(h)}_{L(\;^xg)\cap L(h)})
\eeq
also coincides with the one given in Lemma \ref{indc}.
%}

This shows that \beqn
(H\ot_{L(g)}M)\ot (H\ot_{L(h)}N)\cong \oplus_{x \in \cd}(H \ot_{L(\;^xg\;h)}P(x))\\
\eeqn
as $D_f(U,H)$-modules.

\vskip 0.3 cm
\subsubsection{}
Let $\cb_g$ be the full abelian subcategory of $D_f(U,H)$-mod generated by the set of all irreducible $D_f(U,H)$-modules of the type $H\ot_{L(g)}M$ with $M \in \mi_g$. Recall that for any $M \in \mi_g$ the $D_f(U,H)$-module $H\ot_{L(g)}M$ from Lemma \ref{indc} is denoted by $S_{g,M}$. Note that by Lemma \ref{orbit} one has that $\cb_{\;^xg}=\cb_{g}$ for any $x \in F$.

\md
Fix a set $\Gm=\{g_1,\cdots, g_s\}$ a representatives for the orbits of the action of $F$ on $G$.
Theorem \ref{moddescr} gives a decomposition  
\beq\lb{cdg}
\cg(D_f(U,H)):=\oplus_{g \in \Gm}\cg(\cb_g)
\eeq 
of the Grothendieck ring of $D_f(U,H)$ such that for all $g,h \in G$, one has 
\beq
\cg(\cb_g)\cg(\cb_h)\subseteq \oplus_{x \in D}\cg(\cb_{\;^xgh})
\eeq
where $D$ is a set of double coset representatives of $F_g\backslash F/F_h$.
\subsubsection{} Suppose $f:U \ra H^{*\;\cop}$ is a surjective morphism of Hopf algebras and consider $f^*:F\ra G$ as in Equation \ref{grm}. Since $f^*$ is an injective morphism of groups one can identify the group algebra $\Z F$ as a subalgebra of $\Z G$. Then a basis for $\mtr{C}_{\Z G}(\Z F) $ is given by the orbits sums:
\beq
s(g):=\sum_{x \in F/F_g}f^*(x)gf^*(x)^{-1}
\eeq
Note that in $\Z G$ the product of orbit sums satisfy
\beq\lb{morb}
s(g)s(h)=\sum_{x \in \cd}[F_{\;^xgh}:F_{\;^xg}\cap F_{h}]s(\;^xgh)
\eeq
This formula appears in \cite{scoh}.
\bc
There is a ring surjection
\beq
K_0(D_f(U,H))\ra \mtr{C}_{\Z G}(\Z F) 
\eeq
given by 
\beq 
[S_{g,M}]\mapsto (\dim M) s(g)
\eeq
where $\mtr{C}_{\Z G}(\Z F) $ is the centralizer of $\Z F$ inside the group algebra $\Z G$.
\ec
\bpf
Using Formulae \ref{tep} and \ref{morb} it is easy to check that the above defined map is a ring morphism. Since  $\mtr{C}_{\Z G}(\Z F) $  is $\Z$-spanned by the orbit sums $\co(g)$ it follows that the above morphism is surjective.
\epf
In particular we obtain the following generalization of Corollary 3.2 from \cite{cibils}.
\bc 
For any semisimple Hopf algebra $A$  one has a surjective ring morphism
\beq
K_0(D(A))\ra \mtc{Z}(\Z G)
\eeq
where $\mtc{Z}(\Z G)$ is the center of the group algebra $\Z G$.
\ec
\subsection{On the Grothendieck ring structure of $\cg(D_f(U,H))$ }
In this subsection we show that the Grothendieck ring $\cg(D_f(U,H))$ has the same type ring structure as the rings considered in \cite{scoh} and \cite{bouc}. 
Using the notations from Section 2 of \cite{scoh} one can take $A(g):=K_0(\mtc{B}_g)$ for all $g \in G$ and consider the structure maps 
\beq
m_{g,h}:K_0(\cb_g)\times K_0(\cb_h)\ra K_0(\cb_{gh})
\eeq
 and 
 \beq
c_{g,x}: K_0(\cb_g)\ra K_0(\cb_{\;^xg})
\eeq
induced from the maps defined in Equations \ref{mgh} and respectively Equation \ref{cxg}. Then clearly by Equation \ref{cdg} one has that the invariant ring satisfies $A^L\cong K_0(D_f(U,H))$ as $\Z$-modules since $c_{x,g}$ is the identity map if $x\in F_g$. On the other hand Formula \ref{tep} from the tensor product shows that $A^L$ has the ring structure described in Equation 2.3 from \cite{scoh}. Indeed, in order to realize this fact one has to identify as in page 19 of \cite{scoh}, the module $P(x) \in \ci_{\;^xg_ig_j}$ with the module $\;^yP(x) \in \ci_{g_k}$ where $y\in L$ is chosen such that $\;^{yx}g_i\;^yg_j=g_k$ for a uniquely determined $k=k(x)$. Here $\Gm=\{g_i\}_{i=1}^s\subset G$ is a fixed set of representative elements for the orbits of the action of $F$ on $G$ as above. Then one can use Corollary 2.5 of \cite{scoh}.

%\section{Acknowledgements}This work was supported by the strategic grant POSDRU/89/1.5/S/58852, Project "Postdoctoral programme for training scientific researchers" cofinanced by the European Social Found within the Sectorial Operational Program Human Resources Development 2007 - 2013.
\bibliographystyle{amsplain}
\bibliography{bob}
\end{document}